%% file: main.tex
\documentclass[11pt,a4paper]{article}

\usepackage[a4paper,left=26mm,right=26mm,top=24mm,bottom=24mm,headheight=14pt]{geometry}
\usepackage{amsmath,amssymb,amsthm,mathtools}
\usepackage{newtxmath}
\usepackage{fontspec}
\usepackage{microtype}
\usepackage{fancyhdr}
\usepackage{needspace}
\usepackage{booktabs}
\usepackage{graphicx}
\usepackage{tikz}
\usetikzlibrary{arrows.meta}
\usepackage[font=small,labelfont=bf,labelsep=period]{caption}
\usepackage{placeins}
\usepackage{float}
\usepackage[hidelinks,unicode]{hyperref}
\numberwithin{equation}{section}
\allowdisplaybreaks[1]
\author{%
Xijun Hu\thanks{Email: \texttt{xjhu@sdu.edu.cn}.}\quad
Li Wu\thanks{Email: \texttt{vvvli@sdu.edu.cn}.}\quad
Ao Zhu\thanks{Email: \texttt{zhuao@mail.sdu.edu.cn}.}\\[6pt]
\small School of Mathematics, Shandong University\\
\small Jinan, Shandong 250100, The People's Republic of China
}
\date{}
\newtheorem{theorem}{Theorem}[section]
\newtheorem{lemma}[theorem]{Lemma}
\theoremstyle{remark}
\newtheorem{remark}[theorem]{Remark}
\theoremstyle{plain}
\title{\LARGE Strict Separation of Antiperiodic Degeneracy Curves\\
for Elliptic Euler and Lagrangian Solutions}
\hypersetup{pdftitle={Strict Separation of Antiperiodic Degeneracy Curves for Elliptic Euler and Lagrangian Solutions}}
\begin{document}
\maketitle
\thispagestyle{plain}
\begin{abstract}
We prove strict separation of antiperiodic degeneracy curves associated with elliptic Euler and Lagrangian solutions of the planar three-body problem. For every eccentricity $0<e<1$, each Euler pair issuing from a common circular parameter is strictly ordered in the extended family $\beta>0$, and the two Lagrangian curves remain disjoint. Every antiperiodic degeneracy point in these families with $0<e<1$ has a one-dimensional full kernel. Rotation and reflection decomposition yield, for each model, two half-integer Fourier block operators with identical tails and different first blocks. In the Euler case, a positive definite first-block perturbation gives a definite connecting crossing form; a spectral-flow rectangle and sectorwise index counts then force strict separation. In the Lagrangian case, exact tail elimination and an explicit estimate for a two-dimensional Schur complement exclude simultaneous degeneracy in the two reflection sectors.
\end{abstract}

\section{Introduction and main results}

The study of homographic motions in the three-body problem goes back to the collinear solutions of Euler and the equilateral solutions of Lagrange in the eighteenth century. Both arise from central configurations: the bodies preserve their shape while their common scale and rotation follow a Keplerian motion. For circular Lagrangian solutions, the classical Gascheau--Routh criterion determines the stable mass range. The elliptic problem is more delicate because the linear variational equations have periodic coefficients, and stability depends on both the masses and the eccentricity; see \cite{Roberts2002,MeyerSchmidt2005} for the classical background.

These families also connect the general problem with its restricted limit. When one mass tends to zero while the two primary masses remain positive, the three Euler collinear configurations give $L_1,L_2,L_3$, whereas the two orientations of the Lagrange equilateral configuration give $L_4,L_5$. For $e=0$, these are the five equilibria of the circular restricted problem in the uniformly rotating frame. For $e>0$, they yield the corresponding elliptic homographic motions, which are fixed in rotating-pulsating coordinates. After parameter identification and a constant rotation where necessary, the restricted planar variational equations belong to the same essential operator families considered below. Thus the associated planar linear stability problems are included in this framework \cite{Roberts2002,MeyerSchmidt2005,MSS2006a}.

Modern studies have combined geometric reductions, perturbation methods and numerical analysis. Roberts \cite{Roberts2002} reduced the elliptic Lagrangian stability problem to a four-dimensional system and obtained analytic and numerical information on its stability regions. Meyer and Schmidt \cite{MeyerSchmidt2005} developed symplectic coordinates adapted to central configurations, separating the center-of-mass and Keplerian motions from the essential variational subsystem. Mart\'inez, Sam\`a and Sim\'o investigated the resulting periodic Hamiltonian systems by normal forms near the circular limit \cite{MSS2006a} and by blow-up methods near the singular collision limit \cite{MSS2006b}.

Maslov-type index theory provides a topological approach to these stability questions; see Long's monograph \cite{Long2002} and, for the Maslov index of Lagrangian paths, Robbin and Salamon \cite{RobbinSalamon1993}. For a symplectic fundamental solution $\gamma:[0,2\pi]\to\operatorname{Sp}(4)$ normalized by $\gamma(0)=I_4$ and for $\omega$ on the unit circle, the theory associates an index pair $(i_\omega(\gamma),\nu_\omega(\gamma))$. The nullity $\nu_\omega(\gamma)=\dim_{\mathbb C}\ker(\gamma(2\pi)-\omega I_4)$ detects the corresponding Floquet multiplier, while $i_\omega$ is an intersection index associated with the path and the $\omega$-degeneracy set. Splitting numbers describe the variation of the index with respect to $\omega$, while iteration formulas relate the indices of iterates to the unit-circle spectrum and provide stability criteria. The choice $\omega=-1$ is particularly relevant here because it detects antiperiodic degeneracy.

Variational and symplectic index methods subsequently provided tools for studying the parameter dependence beyond these limiting regimes. Hu and Sun established Morse--Maslov index relations and symmetry-adapted iteration formulas \cite{HuSun2009Sym}, related the variational geometry of central configurations to instability of circular relative equilibria \cite{HuSun2009RE}, and applied Morse indices to elliptic Lagrangian solutions \cite{HuSun2010}. Hu, Long and Sun \cite{HuLongSun2014} studied the Lagrangian family throughout the mass--eccentricity rectangle, identifying its antiperiodic degeneracy curves and the boundary of the hyperbolic region. Hu and Ou \cite{HuOu2016} developed collision-index methods for the near-collision regime. Zhou and Long \cite{ZhouLong2017} established the indexed degeneracy curves and their ordering for the elliptic Euler family.

We briefly recall the origin of the essential equations used here. A planar central configuration $a=(a_1,a_2,a_3)$ generates a homographic solution
\[
q_i(\tau)=r(\theta(\tau))R(\theta(\tau))a_i,
\qquad r(\theta)=\frac{p}{1+e\cos\theta},
\]
where $\tau$ is physical time, $\theta$ is the true anomaly, $p>0$ is the semilatus rectum and $R(\theta)$ is planar rotation. The central-configuration decomposition splits the twelve-dimensional Hamiltonian variational equation into three four-dimensional subsystems: center-of-mass, Kepler and essential \cite{MeyerSchmidt2005}. The essential subsystem governs the nontrivial characteristic multipliers. After normalization in rotating-pulsating coordinates, with the true anomaly denoted by $t$, elimination of the momentum variables gives
\[
-v''-2Jv'+\frac{D_\beta}{1+e\cos t}v=0,
\qquad J=\begin{pmatrix}0&-1\\1&0\end{pmatrix}.
\]
\begin{samepage}
The periodic rotation $z=R(t)v$ gives the self-adjoint equation $A(\beta,e)z=0$, where
\[
A(\beta,e)=-\partial_t^2I_2-I_2+
\frac{R(t)D_\beta R(t)^T}{1+e\cos t},
\qquad
\begin{aligned}
D_\beta^{\mathrm E}&=\operatorname{diag}(2\beta+3,-\beta),\\
D_\beta^{\mathrm L}&=\operatorname{diag}\!\left(
\frac{3+\sqrt{9-\beta}}2,\frac{3-\sqrt{9-\beta}}2\right).
\end{aligned}
\]
\end{samepage}
These are precisely the operators in \cite[Eq.~(2.42)]{ZhouLong2017} and \cite[Eq.~(2.29)]{HuLongSun2014}. Since the coordinate changes are periodic, antiperiodic degeneracy is equivalent to the occurrence of the multiplier $-1$ in the essential monodromy matrix. Its geometric multiplicity equals the dimension of the kernel subject to $z(2\pi)=-z(0)$ and $z'(2\pi)=-z'(0)$.

For these normalized second-order systems, the Morse--Maslov correspondence identifies $i_{-1}(\gamma_{\beta,e})$ with the antiperiodic Morse index of $A(\beta,e)$; see \cite[Eqs.~(2.27) and (2.30)]{HuLongSun2014} and \cite[Eqs.~(2.40) and (2.43)]{ZhouLong2017}. We denote this index by $m^-_a(\beta,e)$ and the kernel dimension by $\nu_a(\beta,e)$, where $a\in\{\mathrm E,\mathrm L\}$ specifies the model. The Morse index counts strictly negative eigenvalues with multiplicity; zero eigenvalues contribute only to the nullity. For these real operators with antiperiodic boundary conditions, the real kernel dimension equals the complex dimension of its complexification.

Throughout, $e\in[0,1)$ is the eccentricity, and $\beta$ follows the convention of the model under consideration. For the Lagrangian problem,
\[
\beta=\frac{27(m_1m_2+m_2m_3+m_3m_1)}{(m_1+m_2+m_3)^2}\in[0,9].
\]
In the triangular restricted problem, this becomes $\beta=27\mu(1-\mu)$, where the normalized primary masses are $1-\mu$ and $\mu$, with $0<\mu<1$.\footnote{At a collinear libration point, let $c=(1-\mu)/r_1^3+\mu/r_2^3$, where $r_1,r_2$ are the distances to the primaries of normalized masses $1-\mu$ and $\mu$, measured in units of their separation. The corresponding diagonal coefficients are $(1+2c,1-c)=(3+2\beta,-\beta)$, so the Euler parameter is $\beta=c-1$.} The Euler parameter is different and has physical range $[0,7]$. Following \cite[Remark~1.2(iii)]{ZhouLong2017}, we consider its mathematical extension to $\beta>0$; the physical interpretation is restricted to curve portions in the physical parameter range. The two parameters are not identified numerically. We suppress model superscripts in the model-specific sections.

For fixed $e$, the monotonicity and jump rules for these Morse indices are established results. Hu--Long--Sun proved that $m^-_{\mathrm L}$ is non-increasing in $\beta$ and decreases by the kernel dimension when a degeneracy is crossed \cite[Lemma~4.4, Corollary~4.5 and Proposition~6.1]{HuLongSun2014}. Zhou--Long proved that $m^-_{\mathrm E}$ is non-decreasing and increases by the kernel dimension \cite[Lemma~4.2 and Corollary~4.3]{ZhouLong2017}. Both indices start at $2$ when $\beta=0$; see \cite[Eq.~(1.7)]{HuLongSun2014} and \cite[Eqs.~(3.4)--(3.5)]{ZhouLong2017}. We use these known index properties to identify and count the branches.

A scalar antecedent of this separation problem arises in the transverse dynamics of the Euler orbit in the spatial isosceles three-body problem. The associated Ince-type Hill operator can be written as
\[
\mathcal S_{\eta,e}=-\frac{d^2}{dt^2}-1-\frac{\eta}{1+e\cos t},
\qquad \eta>0,\quad 0\le e<1,
\]
with antiperiodic boundary conditions, where $\eta$ equals seven times the mass parameter used in~\cite{HuLiuOuQiaoSalomao2026,HuLiuOuSalomaoYu2025} and is distinct from the planar parameters of this paper. For $0<e<1$, each pair of antiperiodic degeneracy curves issuing from a common circular parameter is strictly separated, and the kernel at each degeneracy point is one-dimensional; see~\cite[Theorem~2.9 and Remark~7.16]{HuLiuOuSalomaoYu2025} and~\cite[Proposition~2.1(ii)]{HuLiuOuQiaoSalomao2026}. The latter work further studies eccentricity dependence and proves monotonicity of the transverse rotation number~\cite[Theorem~1.1]{HuLiuOuQiaoSalomao2026}. The present paper continues this investigation of antiperiodic separation for the coupled, two-component essential systems of the planar Euler and Lagrangian families.

For the planar Euler and Lagrangian families, the established curve theory still leaves a distinction between ordering and strict separation. In the Euler problem, the periodic degeneracy curves coincide in pairs and separate successive pairs of antiperiodic curves, but the ordering within each antiperiodic pair is non-strict \cite[Theorem~1.3]{ZhouLong2017}. In the Lagrangian problem, each reflection subspace has a unique degeneracy parameter with a one-dimensional kernel \cite[Section~7]{HuLongSun2014}, but these facts alone allow the two parameters to coincide. The possibility of intersections for $0<e<1$ was explicitly left open in \cite[Remark~1.10]{HuLongSun2014}. We exclude these coincidences in both models.

The main idea is to compare the reflection sectors before locating their degeneracy parameters. In a half-integer Fourier basis, the operators differ only in the lowest-frequency block. For the Euler family, a definite first-block path connects the sectors. A spectral-flow rectangle, with a zero reference on the periodic separators and one simple root per sector in each intervening interval, gives strict ordering without eliminating the common tail. For the Lagrangian family, exact tail elimination and a signed difference of two effective determinants exclude simultaneous degeneracy. Known index counts then identify these roots with the established branches.

Let $\beta_j(-1,e)$ denote the Euler antiperiodic degeneracy parameters ordered increasingly with kernel multiplicity, and let $\beta_s(e)\leq\beta_m(e)$ denote the ordered Lagrangian curves as in \cite{HuLongSun2014}. For a numerical overview, see Figure~\ref{fig:euler-separation} for the Euler degeneracy curves and their periodic separators, and Figure~\ref{fig:lagrange-stability} for the Lagrangian antiperiodic curves and the surrounding stability regions.

\begin{theorem}\label{thm:main}
For every $0<e<1$, the extended Euler essential operator family satisfies
\begin{equation}\label{eq:eulerclaim}
\beta_{2n-1}(-1,e)<\beta_{2n}(-1,e),\qquad n\geq1.
\end{equation}
The elliptic Lagrangian family satisfies
\begin{equation}\label{eq:lagclaim}
\beta_s(e)<\beta_m(e).
\end{equation}
At each of these degeneracy parameters, the full antiperiodic kernel is one-dimensional. At $e=0$, the $n$th Euler pair meets at $\widehat\beta_{n+1/2}$, where
\begin{equation}\label{eq:betahat}
\widehat\beta_q=\frac{q^2-3+\sqrt{9q^4-14q^2+9}}4,\qquad q\geq1,
\end{equation}
whereas the two Lagrangian branches meet at $\beta=3/4$.
\end{theorem}

\begin{remark}[Reflection labels]\label{rem:intro-reflection-labels}
Extend functions antiperiodically to $\mathbb R$ and set
\[
\mathcal T v(t)=\mathcal{N} v(-t),\qquad
\mathcal{N}=\operatorname{diag}(1,-1),\qquad
\mathcal H_\sigma=\ker(\mathcal T-\sigma I),\quad \sigma\in\{+1,-1\}.
\]
In $\mathcal H_+$ the first component is even and the second is odd; $\mathcal H_-$ has the opposite parities. For $0<e<1$, the smaller branches belong to $\mathcal H_+$ in the Euler family and $\mathcal H_-$ in the Lagrangian family, with the larger branches in the opposite sectors. The Euler labels follow from Section~3. For the Lagrangian sector roots, Section~4 recovers the known circular slopes \cite[Proposition~7.4]{HuLongSun2014} in our reflection convention:
\[
\rho^{\mathrm L,\sigma}(e)=\frac34+\sigma\frac{\sqrt{33}}4e+O(e^2),
\]
which fixes the labels near $e=0$; continuity and nonintersection extend them to $0<e<1$. Finally, $\mathcal{N} R(-t)=R(t)\mathcal{N}$ shows that $z=R(t)v$ preserves these labels.
\end{remark}

\begin{remark}[Antiperiodic Morse indices]\label{rem:intro-morse-indices}
Write $m^-_\sigma=m^-(\left.A\right|_{\mathcal H_\sigma})$ and $\nu_\sigma=\dim\ker(\left.A\right|_{\mathcal H_\sigma})$ on the antiperiodic domain. At $\beta=0$, both models have $D_0=\operatorname{diag}(3,0)$ and each sector has $(m^-_\sigma,\nu_\sigma)=(1,0)$ \cite[Eqs.~(7.11) and (7.16)]{HuLongSun2014}. The known index-jump formulas, Theorem~\ref{thm:main} and the preceding reflection labels give, for $0<e<1$ and $n\geq1$,
\begin{center}
\begingroup
\setlength{\tabcolsep}{7pt}
\renewcommand{\arraystretch}{1.22}
\begin{tabular}{llccc}
\toprule
Family & Degeneracy parameter & $m^-$ & $(m^-_+,m^-_-)$ & $(\nu_+,\nu_-)$ \\
\midrule
Euler & $\beta_{2n-1}(-1,e)$ & $2n$ & $(n,n)$ & $(1,0)$ \\
Euler & $\beta_{2n}(-1,e)$ & $2n+1$ & $(n+1,n)$ & $(0,1)$ \\
Lagrangian & $\beta_s(e)$ & $1$ & $(1,0)$ & $(0,1)$ \\
Lagrangian & $\beta_m(e)$ & $0$ & $(0,0)$ & $(1,0)$ \\
\bottomrule
\end{tabular}
\endgroup
\end{center}
Since zero eigenvalues are excluded, the Euler and Lagrangian root values equal their left- and right-hand limits, respectively. At $e=0$, the Euler meeting point $\widehat\beta_{n+1/2}$ has $(m^-,\nu)=(2n,2)$ and sector indices $(n,n)$; the Lagrangian meeting point $\beta=3/4$ has $(m^-,\nu)=(0,2)$ and sector indices $(0,0)$. Both have $(\nu_+,\nu_-)=(1,1)$.
\end{remark}

Section~2 develops the common operator framework; Sections~3 and~4 prove the Euler and Lagrangian results, respectively. The singular endpoint $e=1$ is not included.

\section{A common antiperiodic operator framework}
\subsection{Rotation and reflection decomposition}
Set

\[
R(t)=\begin{pmatrix}\cos t&-\sin t\\\sin t&\cos t\end{pmatrix},\qquad
J=\begin{pmatrix}0&-1\\1&0\end{pmatrix},\qquad
\mathcal{N}=\begin{pmatrix}1&0\\0&-1\end{pmatrix},\qquad S(t)=R(t)\mathcal{N} R(t)^T.
\]

The second-order operators in~\cite[Eq. (2.42)]{ZhouLong2017} and~\cite[Eq. (2.29)]{HuLongSun2014} are, respectively,

\begin{align}
A^{\mathrm E}(\beta,e)&=-\partial_t^2 I_2-I_2+
\frac{(3+\beta)I_2+3(1+\beta)S(t)}{2(1+e\cos t)},\label{eq:AE}\\
A^{\mathrm L}(\beta,e)&=-\partial_t^2 I_2-I_2+
\frac{3I_2+\sqrt{9-\beta}\,S(t)}{2(1+e\cos t)}.\label{eq:AL}
\end{align}

Their common real antiperiodic domain is

\[
\mathcal D_{-1}=\{z\in H^2([0,2\pi];\mathbb R^2):
z(2\pi)=-z(0),\ z'(2\pi)=-z'(0)\}.
\]

The rotated potential matrices are

\begin{equation}\label{eq:D}
\begin{gathered}
D^{\mathrm E}_\beta=\operatorname{diag}(2\beta+3,-\beta),\qquad \beta>0,\\
D^{\mathrm L}_\beta=\operatorname{diag}(a,b),\qquad
a=\frac{3+\sqrt{9-\beta}}2,\quad b=\frac{3-\sqrt{9-\beta}}2,\qquad 0<\beta<9.
\end{gathered}
\end{equation}

We omit the model superscript when working with one model. With $z=Rv$ and $R'=RJ$, we have $z''=R(v''+2Jv'-v)$, and hence

\begin{equation}\label{eq:rotate}
A(\beta,e)(Rv)=RL_{\beta,e}v,\qquad
L_{\beta,e}=-\partial_t^2 I_2-2J\partial_t+\frac{D_\beta}{1+e\cos t}.
\end{equation}

The periodicity of $R$ and $R'$ preserves both antiperiodic boundary conditions. Extend functions antiperiodically to the real line, define $\mathcal Tv(t)=\mathcal{N} v(-t)$, and set $\mathcal H_\sigma=\ker(\mathcal T-\sigma I)$ for $\sigma\in\{+1,-1\}$. Since $J\mathcal{N}=-\mathcal{N} J$ and $D_\beta\mathcal{N}=\mathcal{N} D_\beta$, the operator $L_{\beta,e}$ commutes with $\mathcal T$, so that

\begin{equation}\label{eq:split}
\ker L_{\beta,e}=\ker(L_{\beta,e}|_{\mathcal H_+})
\oplus\ker(L_{\beta,e}|_{\mathcal H_-}).
\end{equation}

For $v=(x,y)^T$, the corresponding boundary conditions on $[0,\pi]$ are

\[
\begin{array}{ll}
\mathcal H_+:&x'(0)=y(0)=x(\pi)=y'(\pi)=0,\\
\mathcal H_-:&x(0)=y'(0)=x'(\pi)=y(\pi)=0.
\end{array}
\]

The identity $\mathcal{N} R(-t)=R(t)\mathcal{N}$ also shows that rotation preserves reflection type.
\subsection{Half-integer Fourier representation and kernel isomorphism}
Set $q_k=k+1/2$ for $k\ge0$. Use the normalized expansions

\begin{align*}
\mathcal H_+:\quad &x=\pi^{-1/2}\sum_{k\ge0}\alpha_k\cos q_kt,
&y=\pi^{-1/2}\sum_{k\ge0}\eta_k\sin q_kt,\\
\mathcal H_-:\quad &x=\pi^{-1/2}\sum_{k\ge0}\alpha_k\sin q_kt,
&y=-\pi^{-1/2}\sum_{k\ge0}\eta_k\cos q_kt.
\end{align*}

The minus sign in the expansion of $y$ on $\mathcal H_-$ makes the differential blocks the same for both reflection types

\begin{equation}\label{eq:C}
C_q=\begin{pmatrix}q^2&2q\\2q&q^2\end{pmatrix},\qquad
C=\operatorname{diag}_{k\ge0}C_{q_k},\qquad
\operatorname{dom}C=\left\{u:\sum_{k\ge0}(1+q_k^4)|u_k|^2<\infty\right\}.
\end{equation}

By Parseval's identity, the normalized Fourier expansions give
\[
\|v\|_{H^2}^2=\sum_{k\ge0}(1+q_k^2+q_k^4)|u_k|^2,
\qquad u_k=(\alpha_k,\eta_k)^T,
\]
where the $H^2$ norm is the sum of the squared $L^2$ norms of $v,v',v''$.
Since
\[
1+q_k^4\le1+q_k^2+q_k^4\le\tfrac32(1+q_k^4),
\]
the Fourier characterization of $H^2$ gives precisely $\operatorname{dom}C$ in~\eqref{eq:C}. Conversely, the corresponding partial sums converge in $H^2$, and continuity of the endpoint traces preserves both antiperiodic boundary conditions.

Let $Q_\sigma$ be the Fourier matrix of multiplication by $1+e\cos t$. The lowest-frequency identities

\[
\cos t\cos(t/2)=\tfrac12\cos(3t/2)+\tfrac12\cos(t/2),\qquad
\cos t\sin(t/2)=\tfrac12\sin(3t/2)-\tfrac12\sin(t/2)
\]

give

\begin{equation}\label{eq:Q}
Q_\sigma=\begin{pmatrix}
I_2+\sigma\frac e2\mathcal{N}&\frac e2 I_2&0&\cdots\\
\frac e2I_2&I_2&\frac e2I_2&\ddots\\
0&\frac e2I_2&I_2&\ddots\\
\vdots&\ddots&\ddots&\ddots
\end{pmatrix},\qquad (1-e)I\le Q_\sigma\le(1+e)I.
\end{equation}

For $k\ge1$, the product-to-sum identities couple $q_k$ only to $q_{k-1}$ and $q_{k+1}$, each with coefficient $e/2$. The displayed lowest-frequency identities account for the term $\sigma(e/2)\mathcal{N}$ in the first block; the sign built into the second basis function of $\mathcal H_-$ does not change the adjacent couplings. The formula first holds for finite sums and extends by density to $\ell^2$. If $v$ corresponds to $u$, then
\[
\langle u,Q_\sigma u\rangle
=\int_0^{2\pi}(1+e\cos t)|v(t)|^2\,dt,
\]
which proves the bounds in~\eqref{eq:Q}.

Thus $Q_\sigma$ is bounded and invertible and commutes with the blockwise action of $D_\beta$. The Fourier representation of the equation is

\begin{equation}\label{eq:Lfourier}
L_\sigma u:=(C+D_\beta Q_\sigma^{-1})u=0,\qquad u\in\operatorname{dom}C.
\end{equation}

The eigenvalues of $C_q$ are $q(q+2)$ and $q(q-2)$. Half-integer frequencies avoid $0$ and $2$, and therefore

\begin{equation}\label{eq:cinverse}
C_{q_k}^{-1}=\begin{pmatrix}c_k&d_k\\d_k&c_k\end{pmatrix},\qquad
c_k=\frac1{q_k^2-4},\qquad d_k=-\frac2{q_k(q_k^2-4)}.
\end{equation}

Since $\|C_{q_k}^{-1}\|=O(q_k^{-2})$, the operator $C^{-1}$ is bounded and self-adjoint on $\ell^2(\mathbb N_0;\mathbb R^2)$ and maps continuously into $\operatorname{dom}C$. Define

\begin{equation}\label{eq:H}
H_\sigma=C^{-1}+D_\beta^{-1}Q_\sigma.
\end{equation}

Since $C^{-1}$ is block diagonal, adding its blocks to those of $D_\beta^{-1}Q_\sigma$ in~\eqref{eq:Q} gives
\[
H_\sigma=\begin{pmatrix}
C_{1/2}^{-1}+D_\beta^{-1}+\sigma\frac e2D_\beta^{-1}\mathcal{N}&\frac e2D_\beta^{-1}&0&\cdots\\
\frac e2D_\beta^{-1}&C_{3/2}^{-1}+D_\beta^{-1}&\frac e2D_\beta^{-1}&\ddots\\
0&\frac e2D_\beta^{-1}&C_{5/2}^{-1}+D_\beta^{-1}&\ddots\\
\vdots&\ddots&\ddots&\ddots
\end{pmatrix}.
\]
Only the first block depends on the reflection sign $\sigma$; all remaining blocks coincide.

The commutation of $D_\beta^{-1}$ and $Q_\sigma$ makes $H_\sigma$ bounded and self-adjoint on all of $\ell^2$. If $u\in\ker L_\sigma$, set $f=Cu$. Then $Q_\sigma f+D_\beta u=0$, or $H_\sigma f=0$. Conversely, if $H_\sigma f=0$, set $u=C^{-1}f\in\operatorname{dom}C$. Then $u=-D_\beta^{-1}Q_\sigma f$, whence $D_\beta Q_\sigma^{-1}u=-f=-Cu$. Consequently,

\begin{equation}\label{eq:kerneliso}
\ker L_\sigma\simeq\ker H_\sigma,\qquad u\longmapsto Cu,\quad f\longmapsto C^{-1}f.
\end{equation}

The isomorphism~\eqref{eq:kerneliso} concerns kernels; it does not identify the Morse indices of $L_\sigma$ and $H_\sigma$. The Euler argument below uses the spectral flow and crossing forms of a normalized auxiliary family; no finite negative inertia is assumed for $H_\sigma$.

The two operators differ only in their first block:

\begin{equation}\label{eq:Hdiff}
H_+-H_-=\operatorname{diag}(eD_\beta^{-1}\mathcal{N},0,0,\ldots).
\end{equation}

The sign of $D_\beta^{-1}\mathcal{N}$ distinguishes the two problems: it is positive definite for the Euler family and indefinite for the Lagrangian family.

\section{Euler solutions: a spectral-flow comparison}
\subsection{A definite first-block path}
In this section, $D_\beta=D^{\mathrm E}_\beta$. Put
\[
S_\beta=|D_\beta|^{1/2}
=\operatorname{diag}(\sqrt{2\beta+3},\sqrt\beta),
\qquad D_\beta=S_\beta\mathcal{N} S_\beta.
\]
The matrices $S_\beta$ and $\mathcal{N}$ act blockwise and commute with $Q_\sigma$. For $\beta>0$, the congruence
\begin{equation}\label{eq:EulerK}
\mathscr K_\sigma=S_\beta H_\sigma S_\beta
=\mathcal{N} Q_\sigma+S_\beta C^{-1}S_\beta
\end{equation}
preserves kernel dimensions. By~\eqref{eq:kerneliso}, its kernel is identified with that of $L_\sigma$ through
\begin{equation}\label{eq:Eulerkernel}
u=C^{-1}S_\beta w,\qquad
w=-\mathcal{N} S_\beta Q_\sigma^{-1}u.
\end{equation}
Let $P_0$ embed $\mathbb R^2$ into the zeroth Fourier block, and connect the two sectors by
\begin{equation}\label{eq:Eulerpath}
\mathscr K(\beta,e,t)
=\frac{\mathscr K_++\mathscr K_-}{2}
+t\frac e2P_0P_0^*,\qquad -1\le t\le1.
\end{equation}
The endpoints are $\mathscr K_-$ and $\mathscr K_+$, respectively. This is a norm-smooth self-adjoint Fredholm family for $\beta>0$, $e<1$: its principal part $\mathcal{N}(Q_++Q_-)/2$ is self-adjoint and invertible by this commutation and~\eqref{eq:Q}, while $C^{-1}$ is compact and the first-block perturbation has finite rank.

The adjacent coupling blocks of~\eqref{eq:Eulerpath} are $(e/2)\mathcal{N}$. For $e>0$, a kernel vector with $w_0=0$ must vanish, by successive block equations. Hence every nonzero kernel vector satisfies
\begin{equation}\label{eq:Eulerconnect}
\langle w,\partial_t\mathscr K\,w\rangle
=\frac e2|w_0|^2>0.
\end{equation}
In the mass direction, work with the Fourier operator $L_\sigma(\beta,e)=C+D_\beta Q_\sigma^{-1}$ on $\operatorname{dom}C$, as in~\eqref{eq:Lfourier}. Since $D_\beta=I_2+(\beta+1)\partial_\beta D_\beta$,
\[
L_\sigma(\beta,e)=L_\sigma(-1,e)+(\beta+1)\partial_\beta L_\sigma,
\qquad \partial_\beta L_\sigma=(\partial_\beta D_\beta)Q_\sigma^{-1},
\]
where the derivative is independent of $\beta$. Here $D_{-1}=I_2$, and $L_\sigma(-1,e)>0$ on $\operatorname{dom}C$ by~\cite[Lemma~4.1(ii)]{ZhouLong2017} and the unitary Fourier representation.

Indeed, let $v\in\mathcal D_{-1}$ be nonzero, set $z(\theta)=R(\theta)v(\theta)$ and $q_e(\theta)=1+e\cos\theta$, and use primes for differentiation in $\theta$. Since $q_e''/q_e=-1+q_e^{-1}$, rotation and integration by parts give
\[
\begin{aligned}
\langle v,L_{-1,e}v\rangle_{L^2}
&=\int_0^{2\pi}\left(|z'|^2+\frac{q_e''}{q_e}|z|^2\right)\,d\theta\\
&=\int_0^{2\pi}q_e^2\left|\left(\frac{z}{q_e}\right)'\right|^2\,d\theta>0.
\end{aligned}
\]
The boundary term vanishes by antiperiodicity. Equality would make $z/q_e$ constant, which is possible only for $z=0$ because $q_e$ is periodic and $z$ is antiperiodic.

For $w\in\ker\mathscr K_\sigma$, let $u$ be given by~\eqref{eq:Eulerkernel}; then $L_\sigma(\beta,e)u=0$. Differentiating~\eqref{eq:EulerK} and using~\eqref{eq:Eulerkernel} together with $\partial_\beta D_\beta=2\mathcal{N} S_\beta\partial_\beta S_\beta$ give
\begin{equation}\label{eq:Eulermass}
\begin{aligned}
\langle w,\partial_\beta\mathscr K_\sigma w\rangle
&=2\langle(\partial_\beta S_\beta)w,u\rangle
=-\langle(\partial_\beta D_\beta)Q_\sigma^{-1}u,u\rangle\\
&=-\langle(\partial_\beta L_\sigma)u,u\rangle
=\frac{\langle L_\sigma(-1,e)u,u\rangle}{\beta+1}>0
\qquad(w\ne0).
\end{aligned}
\end{equation}
The last equality uses the affine identity above on $\ker L_\sigma(\beta,e)$.

We use the convention that spectral flow counts crossings from negative to positive as positive. For its index-theoretic origins, see Atiyah, Patodi and Singer~\cite{AtiyahPatodiSinger1976}. For spectral flow of bounded self-adjoint Fredholm paths, including additivity and homotopy invariance with invertible endpoints, we use~\cite{Phillips1996}; for the classical crossing-form approach, see~\cite[Section~4]{RobbinSalamon1995}. A variational treatment for strongly indefinite functionals is given by Fitzpatrick, Pejsachowicz and Recht~\cite{FitzpatrickPejsachowiczRecht1999}. In the present bounded setting, a local spectral projection onto the finite-dimensional spectral subspace near zero reduces each regular crossing to a matrix path. Its contribution is the signature of the crossing form, so~\eqref{eq:Eulerconnect} and~\eqref{eq:Eulermass} each contribute the kernel dimension.

\subsection{Periodic separators and the sector counts}
Fix $n\ge1$ and $e_*\in(0,1)$. The continuity and ordering results in~\cite[Theorem~1.3, Lemma~4.5, Theorem~4.11 and Corollary~4.12]{ZhouLong2017} allow us to choose
\[
0<\varepsilon<\min_{0\le e\le e_*}\beta_1(-1,e)
\]
and to consider the moving interval $I_n(e)=(a_n(e),b_n(e))$, where
\begin{equation}\label{eq:moving}
a_n(e)=\begin{cases}\varepsilon,&n=1,\\
\beta_{2n-2}(1,e),&n\ge2,
\end{cases}
\qquad b_n(e)=\beta_{2n}(1,e).
\end{equation}
Both sector operators are invertible at its endpoints, by periodic--antiperiodic nonintersection and the choice of $\varepsilon$ at the first left endpoint. The known ordering also places precisely the $n$th antiperiodic pair in this interval:
\begin{equation}\label{eq:knownorder}
a_n(e)<\beta_{2n-1}(-1,e)\le\beta_{2n}(-1,e)<b_n(e).
\end{equation}
The positive choice of $\varepsilon$ keeps all homotopies away from $\beta=0$, where $S_\beta$ is only continuous.

At $e=0$, the Fourier operator is block diagonal, with the same blocks $C_q+D_\beta$ in both reflection sectors. A direct calculation gives
\[
\det(C_q+D_\beta)
=-2\beta^2+(q^2-3)\beta+q^2(q^2-1).
\]
For $q=n+1/2$, its unique positive root is $\widehat\beta_q$ from~\eqref{eq:betahat}, and this root is simple. The circular ordering places it in $I_n(0)$ and excludes roots from all other half-integer frequencies in that interval; see~\cite[Section~3.2]{ZhouLong2017}. Thus each sector has exactly one simple root $\widehat\beta_{n+1/2}$ in $I_n(0)$. Homotopy invariance along these moving intervals and~\eqref{eq:Eulermass} give
\begin{equation}\label{eq:Eulerunit}
\operatorname{sf}_{\beta:a_n(e)\to b_n(e)}\mathscr K_\sigma(\beta,e)=1.
\end{equation}
Since every crossing is positive by~\eqref{eq:Eulermass}, each sector has exactly one simple root in $I_n(e)$, denoted by $\rho_n^\sigma(e)$.

There is also a zero reference for the connecting path at the left endpoint. As $e$ varies from $0$ to $e_*$ along $\beta=a_n(e)$, the two sector endpoints remain invertible. At $e=0$, the connecting path is constant, so
\begin{equation}\label{eq:Eulerwall}
\operatorname{sf}_{t:-1\to1}\mathscr K(a_n(e_*),e_*,t)=0.
\end{equation}

\subsection{Strict separation from the spectral-flow rectangle}
\begin{figure}[H]
\centering
\input{figures/euler-proof-panels.tex}
\caption{Schematic of the Euler separation argument at fixed $0<e<1$. (a) The arrows orient the boundary of $[a,\beta]\times[-1,1]$; edge labels are the corresponding spectral-flow contributions. The left reference is zero, so $\Phi(\beta)=N_+(\beta)-N_-(\beta)\ge0$. (b) Only a neighbourhood of the simple zero at $(\rho^-,-1)$ is shown. The local zero curve moves to $\beta_*<\rho^-$ as $t$ increases. Its interior crossing gives $\Phi(\beta_*)\ge1$, and $N_-(\beta_*)=0$ then forces $\rho^+<\beta_*$. Here $t$ is the auxiliary connecting parameter.}
\label{fig:euler-proof}
\end{figure}
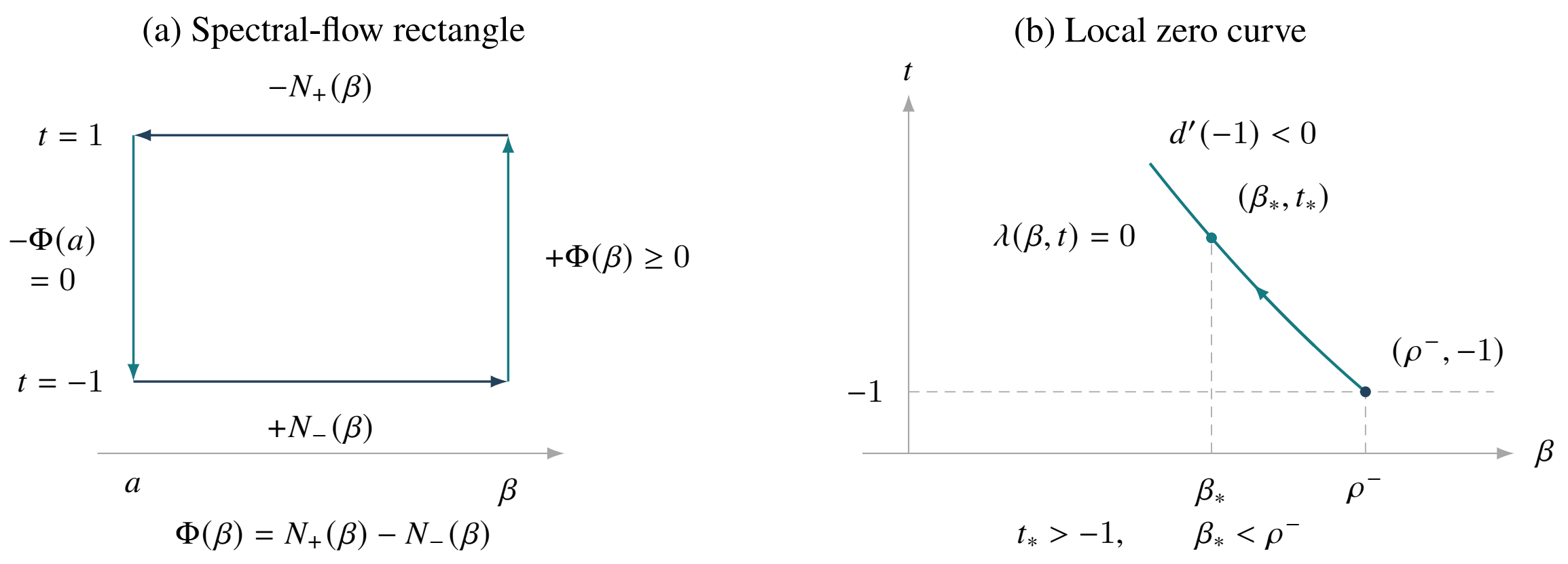

The rectangle count and the local displacement that yields strict separation are illustrated schematically in Figure~\ref{fig:euler-proof}.

Fix $e=e_*$ and abbreviate $a=a_n(e)$, $\rho^\sigma=\rho_n^\sigma(e)$. For $\beta\in I_n(e)$ different from both sector roots, set
\[
N_\sigma(\beta)=\operatorname{sf}_{s:a\to\beta}\mathscr K_\sigma(s,e)
=\mathbf1_{\{\beta>\rho^\sigma\}},\qquad
\Phi(\beta)=\operatorname{sf}_{t:-1\to1}\mathscr K(\beta,e,t).
\]
Traversing the boundary of $[a,\beta]\times[-1,1]$ along the bottom, right, top and left edges, respectively, gives zero total spectral flow, hence
\[
N_-(\beta)+\Phi(\beta)-N_+(\beta)-\Phi(a)=0.
\]
Since $\Phi(a)=0$ by~\eqref{eq:Eulerwall}, this gives
\begin{equation}\label{eq:Eulerrectangle}
\Phi(\beta)=N_+(\beta)-N_-(\beta)\ge0,
\end{equation}
where the inequality follows from~\eqref{eq:Eulerconnect}.

Consider the simple zero of $\mathscr K(\beta,e,-1)=\mathscr K_-(\beta,e)$ at $\beta=\rho^-$. Its local eigenvalue continuation $\lambda(\beta,t)$ satisfies $\partial_\beta\lambda>0$ and $\partial_t\lambda>0$ at $(\rho^-,-1)$ by~\eqref{eq:Eulermass} and~\eqref{eq:Eulerconnect}. The implicit function theorem therefore gives a local zero curve $\beta=d(t)$ with
\[
d(-1)=\rho^-,\qquad
d'(-1)=-\frac{\partial_t\lambda}{\partial_\beta\lambda}<0.
\]
Since $d$ is strictly decreasing near $-1$, choose $t_*\in(-1,1)$ sufficiently close to $-1$ so that $\beta_*=d(t_*)\in I_n(e)$, $\beta_*<\rho^-$, and $\beta_*\ne\rho^+$. Both sector endpoints at $\beta_*$ are then invertible, while the connecting path has an interior crossing. All its crossings are positive, and $N_-(\beta_*)=0$, so~\eqref{eq:Eulerrectangle} yields
\[
1\le\Phi(\beta_*)=N_+(\beta_*)\le1.
\]
Thus $N_+(\beta_*)=1$, and hence $\rho^+<\beta_*<\rho^-$.

Consequently, the kernel correspondence and~\eqref{eq:knownorder} identify
\[
\beta_{2n-1}(-1,e)=\rho_n^+(e)
<\rho_n^-(e)=\beta_{2n}(-1,e).
\]
Each sector root is simple and the two sectors never vanish together, so the full antiperiodic kernel is one-dimensional at each root. Since $n$ and $e_*<1$ were arbitrary, this proves the Euler part of Theorem~\ref{thm:main}.

\begin{figure}[H]
\centering
\includegraphics[width=0.92\textwidth]{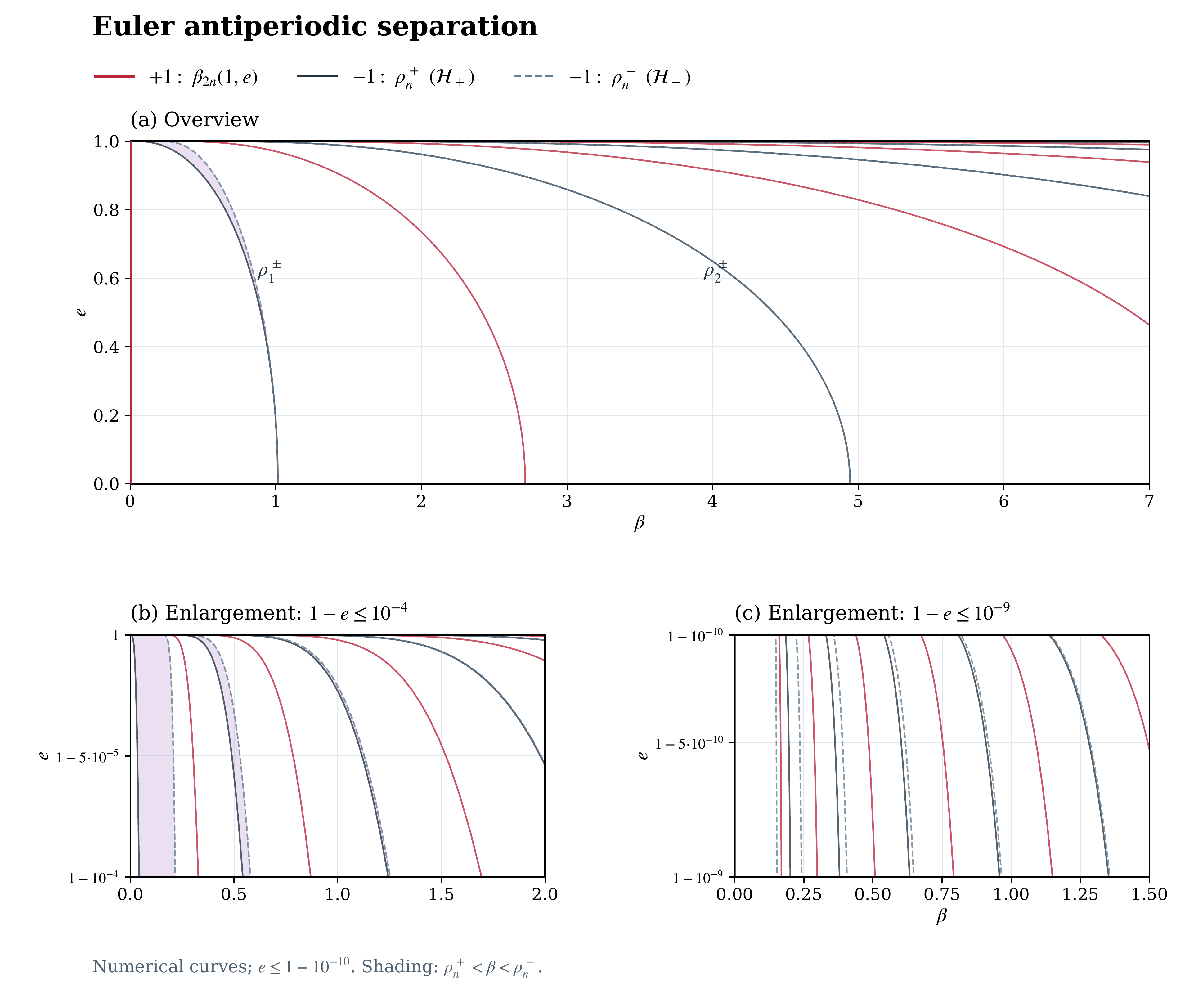}
\caption{Numerical Euler degeneracy curves in the physical range $0\le\beta\le7$. Red curves are the periodic separators $\beta_{2n}(1,e)$; the red boundary $\beta=0$ is also periodic-degenerate. Solid dark and dashed grey curves represent $\rho_n^+=\beta_{2n-1}(-1,e)$ in $\mathcal H_+$ and $\rho_n^-=\beta_{2n}(-1,e)$ in $\mathcal H_-$, respectively. Shading marks the intervals between these paired branches. Panels (b) and (c) enlarge the near-$1$ eccentricity range, with linear axes throughout. Curves interpolate computed roots up to $e=1-10^{-10}$; gaps below numerical or graphical resolution may appear closed. No endpoint at $e=1$ is asserted.}
\label{fig:euler-separation}
\end{figure}

Figure~\ref{fig:euler-separation} illustrates the ordering just proved, together with the periodic separators used in~\eqref{eq:moving}. The enlargements display the same numerical branches near $e=1$; the argument above establishes separation independently of their graphical resolution. The plotted roots were computed by shooting the coupled variational equations on a half period, using an exterior-square formulation and bracketed root searches. Near $e=1$, an exact change of independent variable resolves the rapidly varying coefficients. Refinement of the integration and root-search tolerances, together with an independent QR-based shooting check at selected roots, provides numerical consistency checks; these computations are floating-point approximations.

\FloatBarrier

\section{Lagrangian solutions: scalar control of the common tail}
\subsection{The effective first block and the separation criterion}
Fix $0<\beta<9$ and $0<e<1$, and take $D_\beta=D^{\mathrm L}_\beta=\operatorname{diag}(a,b)$. Thus
$a+b=3$, $a-b=\sqrt{9-\beta}$, and $ab=\beta/4$.
The congruence
\begin{equation}\label{eq:K}
K_\sigma=D_\beta^{1/2}H_\sigma D_\beta^{1/2}
=Q_\sigma+D_\beta^{1/2}C^{-1}D_\beta^{1/2}
\end{equation}
preserves kernel dimensions and makes every adjacent coupling equal to $(e/2)I_2$. Its common diagonal blocks are
\begin{equation}\label{eq:B}
B_k=I_2+D_\beta^{1/2}C_{q_k}^{-1}D_\beta^{1/2}
=\begin{pmatrix}
1+ac_k&\frac{\sqrt\beta}{2}d_k\\
\frac{\sqrt\beta}{2}d_k&1+bc_k
\end{pmatrix},
\end{equation}
so that
\[
K_\sigma=\begin{pmatrix}
B_0+\sigma\frac e2\mathcal{N}&\frac e2I_2&0&\cdots\\
\frac e2I_2&B_1&\frac e2I_2&\ddots\\
0&\frac e2I_2&B_2&\ddots\\
\vdots&\ddots&\ddots&\ddots
\end{pmatrix}.
\]
Since the first-block difference $e\mathcal N$ is indefinite, the positivity-based spectral-flow argument used in the Euler case does not directly apply. Instead, we eliminate the common tail and compare the determinants of the resulting two-dimensional Schur complements to exclude simultaneous degeneracy in the two reflection sectors.

Let $T_k$ be the common tail compression to blocks $k,k+1,\ldots$, and let $E_k$ embed $\mathbb R^2$ into its first block. These compressions are bounded and self-adjoint. The estimate needed for separation is stated first; its proof occupies the next two subsections.
\begin{lemma}\label{lem:effective}
For $0<\beta<9$ and $0<e<1$, the common tail $T_1$ is boundedly invertible. Define
\begin{equation}\label{eq:G}
G_0=B_0-\frac{e^2}{4}E_1^*T_1^{-1}E_1,
\qquad G_\sigma=G_0+\sigma\frac e2\mathcal{N}.
\end{equation}
Then $\ker K_\sigma\simeq\ker G_\sigma$, and
\begin{equation}\label{eq:effectivegap}
(G_0)_{11}-(G_0)_{22}
\le-\frac4{15}(a-b)(1-e^2)<0.
\end{equation}
\end{lemma}

Indeed, once this lemma is proved, the identity
\begin{equation}\label{eq:detdiff}
\det G_+-\det G_-=e\big((G_0)_{22}-(G_0)_{11}\big)
\end{equation}
excludes simultaneous degeneracy in the two sectors. The estimate has a specific purpose: before elimination,
\[
(B_0)_{11}-(B_0)_{22}=-\frac4{15}(a-b).
\]
We must show that the common tail cannot cancel this negative diagonal difference. The proof bounds the increase caused by the tail by $4e^2(a-b)/15$, leaving the strict margin in~\eqref{eq:effectivegap}. To do this, we first control a normalized coefficient of the positive high-frequency tail by a scalar recursion, and then transfer that control through blocks $1$ and $0$.

\subsection{Scalar comparison for the positive high-frequency tail}
For $k\ge2$, $q_k>2$ and $C_{q_k}^{-1}>0$. The compression bounds in~\eqref{eq:Q} give
\begin{equation}\label{eq:Tpositive}
T_k=Q_{\ge k}+\operatorname{diag}_{j\ge k}
\big(D_\beta^{1/2}C_{q_j}^{-1}D_\beta^{1/2}\big)
\ge(1-e)I>0.
\end{equation}
The same bound holds for every finite tail compression. Thus all these positive tails are boundedly invertible, with inverse norm at most $(1-e)^{-1}$.

Retain blocks $1$ through $N$, where $N\ge2$. Successive Schur complements give the effective blocks
\begin{equation}\label{eq:Mrec}
M_N^{(N)}=B_N,\qquad
M_k^{(N)}=B_k-\frac{e^2}{4}\big(M_{k+1}^{(N)}\big)^{-1},
\quad k=N-1,\ldots,1.
\end{equation}
Every inverse here concerns a positive tail starting at $k+1\ge2$. No inverse of $M_1^{(N)}$ is required at this stage.

Suppress the superscript temporarily and write
\begin{equation}\label{eq:Mform}
M_k=\begin{pmatrix}
U_k+aV_k&\frac{\sqrt\beta}{2}W_k\\
\frac{\sqrt\beta}{2}W_k&U_k+bV_k
\end{pmatrix},\qquad
\Delta_k=\det M_k=U_k(U_k+3V_k)+\frac\beta4(V_k^2-W_k^2).
\end{equation}
Here $V_k$ records the diagonal difference, since $(M_k)_{11}-(M_k)_{22}=(a-b)V_k$. The other two coefficients will provide the signs and denominator bounds needed to control this difference. The two-dimensional inverse formula and $a+b=3$ yield
\begin{equation}\label{eq:UVW}
\begin{aligned}
U_k&=1-\frac{e^2(U_{k+1}+3V_{k+1})}{4\Delta_{k+1}},\\
V_k&=c_k+\frac{e^2V_{k+1}}{4\Delta_{k+1}},\qquad
W_k=d_k+\frac{e^2W_{k+1}}{4\Delta_{k+1}},
\end{aligned}
\end{equation}
with terminal values $U_N=1$, $V_N=c_N$, and $W_N=d_N$.

Only a uniform lower bound of $1/2$ is needed. Backward induction gives, for $2\le k\le N$,
\begin{equation}\label{eq:highsigns}
\frac12\le U_k\le1,\qquad V_k>0,\qquad -V_k<W_k<0,
\qquad \Delta_k\ge U_k(U_k+3V_k)>0.
\end{equation}
For $k\ge2$, we have $c_k>0$, $d_k<0$, and $c_k+d_k=[q_k(q_k+2)]^{-1}>0$. Thus the statements hold at $N$. If they hold at $k+1$, then
\[
U_k\ge1-\frac{e^2}{4U_{k+1}}\ge1-\frac{e^2}{2}>\frac12,
\qquad
V_k+W_k=c_k+d_k+\frac{e^2(V_{k+1}+W_{k+1})}{4\Delta_{k+1}}>0.
\]
The other signs follow directly from~\eqref{eq:UVW}. In particular, $V_k^2-W_k^2>0$ makes the last term in the determinant formula~\eqref{eq:Mform} nonnegative.

The quantity passed from the high-frequency tail to block $1$ is the ratio
\begin{equation}\label{eq:theta}
\theta_k=\frac{V_k}{U_k+3V_k},\qquad 0<\theta_k<\frac13,\quad k\ge2.
\end{equation}
Indeed, eliminating the common factor in the first two recursions gives
\begin{equation}\label{eq:theta-transfer}
V_k=c_k+(1-U_k)\theta_{k+1}.
\end{equation}
Thus a bound on $\theta_2$ will directly constrain the relation between $U_1$ and $V_1$. We now obtain this bound without estimating $V_2$ separately.

An explicitly solvable boundary recurrence suggests the comparison function. At the coefficient values $\beta=0$, $e=1$, the diagonal recurrence has the solution $U_*=1/2$ and
\[
X_*(q)=1+\frac3{q^2-4}-\frac1{4X_*(q+1)}
=\frac{q^2-1}{2q(q-2)},\qquad q>2.
\]
Writing $X_*=U_*+3V_*$ gives the normalized coefficient $V_*/X_*$, which is the function $\Theta$ in the following lemma. We use this explicit sequence to control the tail for $e<1$; no bounded tail inverse at $e=1$ is assumed.

\begin{lemma}\label{lem:high}
For every finite tail and $2\le k\le N$,
\begin{equation}\label{eq:thetabound}
0<\theta_k\le\Theta(q_k),\qquad
\Theta(q)=\frac{2q-1}{3(q^2-1)},\quad q>2.
\end{equation}
In particular, $\theta_2\le16/63$.
\end{lemma}
\begin{proof}
For $k<N$, equations~\eqref{eq:theta-transfer} and~\eqref{eq:highsigns} imply
\[
\frac{V_k}{U_k}
=\frac{c_k+\theta_{k+1}}{U_k}-\theta_{k+1}
\le2(c_k+\theta_{k+1})-\theta_{k+1}
=2c_k+\theta_{k+1}.
\]
Since $z\mapsto z/(1+3z)$ is increasing for $z\ge0$,
\begin{equation}\label{eq:scalarcomparison}
\theta_k\le\Phi_{q_k}(\theta_{k+1}),\qquad
\Phi_q(t)=\frac{2c(q)+t}{1+6c(q)+3t},\quad c(q)=\frac1{q^2-4}.
\end{equation}
The map $\Phi_q$ is increasing, and direct substitution gives
\begin{equation}\label{eq:Thetarec}
\Theta(q)=\Phi_q\big(\Theta(q+1)\big).
\end{equation}
At the terminal block,
\[
\theta_N=\frac{c_N}{1+3c_N}
=\frac1{q_N^2-1}<\Theta(q_N),
\]
because $q_N>2$. Backward induction using the monotonicity of $\Phi_q$ and~\eqref{eq:Thetarec} proves~\eqref{eq:thetabound}. Finally, $\Theta(5/2)=16/63$.
\end{proof}

These bounds pass to the infinite tail by a short Galerkin argument. Fix $k\ge2$, let $P_N$ project onto blocks $k,\ldots,N$, and let $T_{k,N}$ be the corresponding compression. For $f$ supported in block $k$, put $y=T_k^{-1}f$ and $y_N=T_{k,N}^{-1}P_Nf$, extended by zero. The two equations imply Galerkin orthogonality:
\[
\langle T_k(y-y_N),z\rangle=0
\qquad (z\in\operatorname{Ran}P_N).
\]
Thus $y_N$ is the best approximation to $y$ from $\operatorname{Ran}P_N$ in the energy norm $\|z\|_{T_k}^2=\langle T_kz,z\rangle$. Comparing with $P_Ny$ and using~\eqref{eq:Tpositive}, we obtain
\[
\|y-y_N\|^2
\le\frac{\|y-y_N\|_{T_k}^2}{1-e}
\le\frac{\|T_k\|}{1-e}\|y-P_Ny\|^2\longrightarrow0.
\]
Thus the first inverse blocks converge, and for each fixed $k\ge1$,
\begin{equation}\label{eq:exactM}
M_k^{(N)}\longrightarrow M_k
:=B_k-\frac{e^2}{4}E_{k+1}^*T_{k+1}^{-1}E_{k+1}.
\end{equation}
For $k\ge2$, minimizing the tail quadratic form over the remaining blocks shows that the Schur complements satisfy $M_k^{(N)},M_k\ge(1-e)I_2$. Their inverses therefore converge as well, and block inversion gives $E_k^*T_k^{-1}E_k=M_k^{-1}$. The coefficient representation, recursions, and non-strict bounds therefore pass to the limit. The strict signs $V_k>0$, $W_k<0$, and $V_k+W_k>0$ follow again from~\eqref{eq:UVW}. In particular,
\begin{equation}\label{eq:theta2}
U_2\ge\frac12,\qquad 0<\theta_2\le\frac{16}{63}.
\end{equation}
The recursions also hold at block $1$, where no invertibility has yet been assumed.

\subsection{Transfer through the two low-frequency blocks}
\begin{proof}[Proof of Lemma~\ref{lem:effective}]
At block $1$, $c_1=-4/7$ and $d_1=16/21$. Write
\[
x=U_1=1-\frac{e^2(U_2+3V_2)}{4\Delta_2},\qquad y=-V_1.
\]
The high-frequency determinant bound gives
\[
\frac12<1-\frac{e^2}{2}\le1-\frac{e^2}{4U_2}\le x<1.
\]
Equation~\eqref{eq:theta-transfer}, now at $k=1$, and~\eqref{eq:theta2} give the useful linear constraint
\begin{equation}\label{eq:lowlinear}
y=\frac47-\theta_2(1-x)
\ge\frac47-\frac{16}{63}(1-x)=\frac{20+16x}{63}.
\end{equation}
In particular, $y>0$ and $3y-x\ge\frac{15(4-x)}{63}>0$. Also,
\[
W_1-y=W_1+V_1
=\frac4{21}+\frac{e^2(W_2+V_2)}{4\Delta_2}>0.
\]
Hence the determinant at block $1$ satisfies
\begin{equation}\label{eq:delta1}
-\Delta_1=x(3y-x)+\frac\beta4(W_1^2-y^2)
\ge x(3y-x)>0.
\end{equation}
This proves invertibility of $M_1$. Since $T_2$ is positive and boundedly invertible, block elimination gives a bounded invertible congruence between $T_1$ and $M_1\oplus T_2$. Therefore $T_1$ is boundedly invertible, and
\[
E_1^*T_1^{-1}E_1=M_1^{-1}.
\]

The same linear constraint controls the correction passed to block $0$. For fixed $x$, the ratio $\frac{y}{x(3y-x)}$ decreases with $y$ on $3y>x$. Substituting the lower bound for $y$ from~\eqref{eq:lowlinear} and using~\eqref{eq:delta1} give
\begin{equation}\label{eq:ratio}
0<\frac{V_1}{\Delta_1}
=\frac{y}{-\Delta_1}
\le\frac{y}{x(3y-x)}
\le\frac{20+16x}{15x(4-x)}\le\frac{16}{15}.
\end{equation}
The last inequality is equivalent to $(2x-1)(5-2x)\ge0$, which holds for $1/2<x<1$.

We can now eliminate the entire tail of $K_\sigma$:
\[
K_\sigma=\begin{pmatrix}
B_0+\sigma\frac e2\mathcal{N}&\frac e2E_1^*\\
\frac e2E_1&T_1
\end{pmatrix},\qquad
G_0=B_0-\frac{e^2}{4}M_1^{-1}.
\]
The tail equation uniquely determines $\zeta=-(e/2)T_1^{-1}E_1\xi$, giving the kernel isomorphism
\[
\xi\longmapsto\big(\xi,-(e/2)T_1^{-1}E_1\xi\big),
\qquad \ker G_\sigma\simeq\ker K_\sigma.
\]
The two-dimensional inverse formula applied to~\eqref{eq:Mform} gives the same coefficient recursion at block $0$. Since $c_0=-4/15$, the diagonal difference of $G_0$ is $(a-b)V_0$, with
\begin{equation}\label{eq:V0}
V_0=-\frac4{15}+\frac{e^2V_1}{4\Delta_1}
\le-\frac4{15}(1-e^2)<0.
\end{equation}
This proves~\eqref{eq:effectivegap}. The bound shows explicitly that the common tail increases the diagonal-difference coefficient by at most $4e^2/15$, so it cannot cancel the initial value $-4/15$ for $e<1$.
\end{proof}

\subsection{Identification and ordering of the degeneracy curves}

\begin{figure}[H]
\centering
\includegraphics[width=\textwidth]{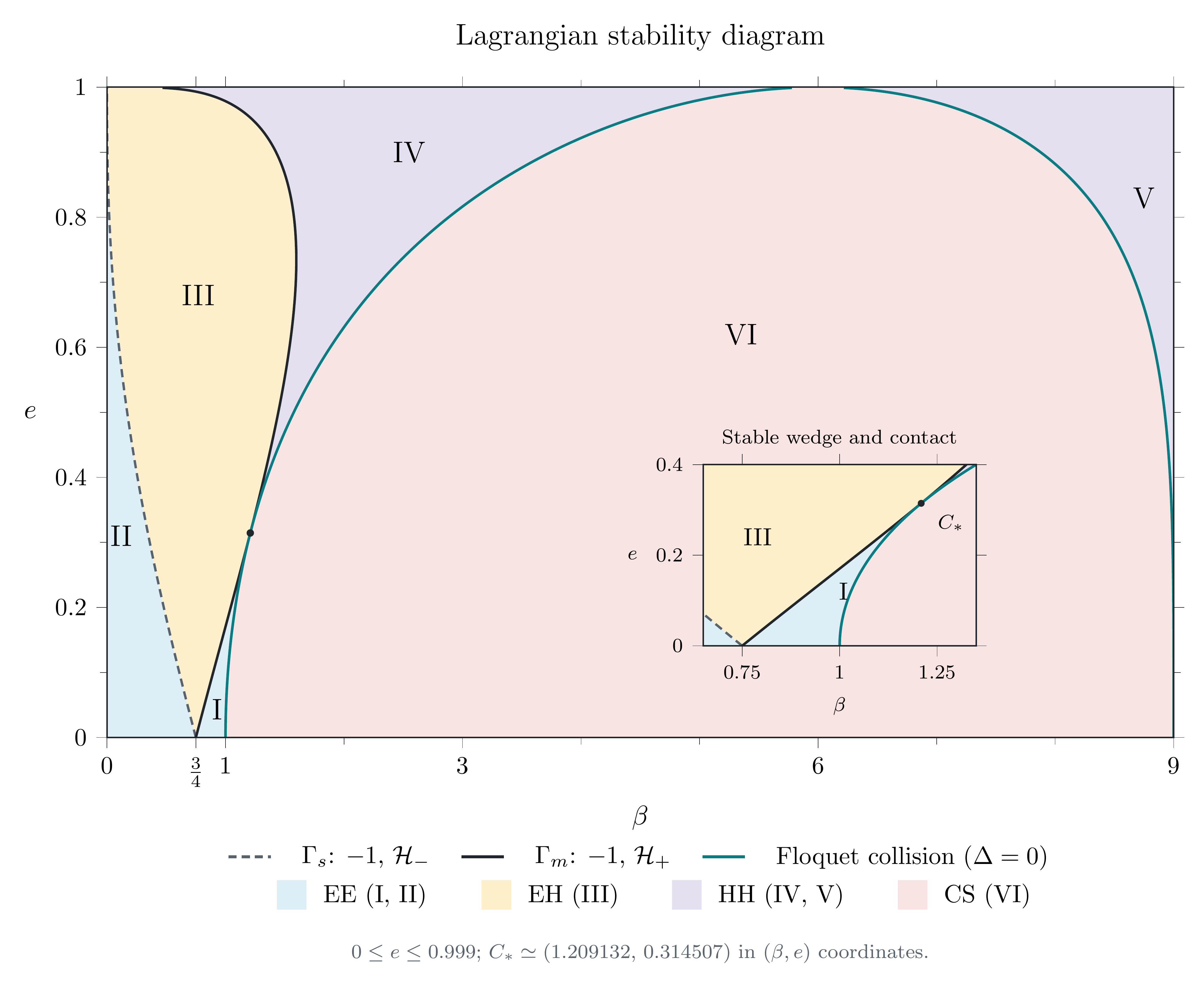}
\caption{Numerical Lagrangian stability diagram for $0\le\beta\le9$ and $0\le e\le0.999$. The dashed curve $\Gamma_s$ is $\beta=\beta_s(e)$ in $\mathcal H_-$, and the solid dark curve $\Gamma_m$ is $\beta=\beta_m(e)$ in $\mathcal H_+$. They meet at $(3/4,0)$ and are strictly separated for $0<e<1$. The teal curve indicates the Floquet collision condition $\Delta_{\mathrm F}=0$, denoted by $\Delta=0$ in the plot. The marked contact $C_*\simeq(1.209132,0.314507)$ is numerical and is not an intersection of $\Gamma_s$ and $\Gamma_m$. Region labels describe the essential Floquet spectrum: EE, elliptic--elliptic; EH, elliptic--hyperbolic; HH, hyperbolic--hyperbolic; CS, complex saddle. The background classification and the displayed contact coordinates are numerical information; the digits shown are not certified error bounds.}
\label{fig:lagrange-stability}
\end{figure}

Figure~\ref{fig:lagrange-stability} places the two separated Lagrangian curves in a numerical stability diagram. The inset distinguishes their common circular starting point from the contact of the upper antiperiodic curve with a Floquet collision curve; the latter is a different spectral condition. For the essential monodromy matrix $M=\gamma_{\beta,e}(2\pi)$, write
\[
\det(\lambda I_4-M)=\lambda^4-\tau\lambda^3+\kappa\lambda^2-\tau\lambda+1.
\]
The substitution $z=\lambda+\lambda^{-1}$ gives $z^2-\tau z+\kappa-2=0$. Thus the Floquet collision discriminant is $\Delta_{\mathrm F}=\tau^2-4(\kappa-2)$; the symbol $\Delta$ in the figure denotes this discriminant. The numerical region classification and contact have antecedents in~\cite{MSS2006a,MSS2006b}; see also the discussion in~\cite[Section~10]{HuLongSun2014}.

Lemma~\ref{lem:effective} and~\eqref{eq:detdiff} give

\begin{equation}\label{eq:detpositive}
\det G_+-\det G_-
\ge\frac{4e}{15}\sqrt{9-\beta}\,(1-e^2)>0.
\end{equation}

Thus the two reflection kernels cannot be nontrivial at the same parameter.

To identify the curves, use the root count in each reflection subspace from~\cite[Corollary~4.3, Lemma~4.4, Proposition~7.1 and Eqs.~(7.7), (7.11), (7.16)]{HuLongSun2014}. More precisely, $A^{\mathrm L}(9,e)$ is positive definite on the antiperiodic space, and

\[
\frac\partial{\partial\beta}\left(\frac{A^{\mathrm L}(\beta,e)}{\sqrt{9-\beta}}\right)
=\frac{A^{\mathrm L}(9,e)}{2(9-\beta)^{3/2}}>0.
\]

The displayed derivative remains positive definite after restriction to either reflection subspace. Positive rescaling preserves the Morse index and kernel, so every zero crossing in a sector is regular and decreases its Morse index by the kernel dimension. At $\beta=0$, each sector has trivial kernel and exactly one negative direction. Since $A^{\mathrm L}(9,e)>0$, continuity gives positivity for $\beta<9$ sufficiently close to $9$. The total index decrease in each sector is therefore one, using the rescaled family only for $\beta<9$, and all crossings have the same direction. Hence each sector has exactly one degeneracy parameter in $0<\beta<9$, with a one-dimensional kernel. Denote these parameters by $\rho^+(e)$ and $\rho^-(e)$. Equation~\eqref{eq:detpositive} and the kernel isomorphisms give $\rho^+(e)\ne\rho^-(e)$. Since

\[
\beta_s(e)=\min\{\rho^+(e),\rho^-(e)\},\qquad
\beta_m(e)=\max\{\rho^+(e),\rho^-(e)\},
\]

we obtain $\beta_s(e)<\beta_m(e)$, with a one-dimensional full antiperiodic kernel at each curve. At $e=0$, the identity

\[
\det(C_q+D_\beta)=q^4-q^2+\frac\beta4
\]

shows that only $q=1/2$ yields a degeneracy parameter in $0<\beta<9$, namely $\beta=3/4$.

To identify the reflection labels near this circular point, we recover the known slopes from~\cite[Proposition~7.4]{HuLongSun2014} directly from the effective block. At $(\beta,e)=(3/4,0)$, the common tail $T_1$ is block diagonal: $\det B_1=-5/7-4\beta/63<0$, while $B_k\ge I_2$ for $k\ge2$. Thus $T_1$ is boundedly invertible there. Its real-analytic dependence on $(\beta,e)$ extends its invertibility and the real analyticity of its inverse to a neighbourhood of $(3/4,0)$, allowing either sign of $e$. Equation~\eqref{eq:G} consequently gives $G_0=B_0+O(e^2)$, uniformly for $\beta$ near $3/4$. Using $c_0=-4/15$, $d_0=16/15$, and $a-b=\sqrt{9-\beta}$, we obtain
\[
\det G_\sigma(\beta,e)
=\frac{3-4\beta}{15}
+\sigma\frac{2e}{15}\sqrt{9-\beta}+O(e^2).
\]
Since $\partial_\beta\det G_\sigma(3/4,0)=-4/15$, the real-analytic implicit function theorem gives a local branch agreeing with $\rho^\sigma(e)$ for small positive $e$, with
\[
\rho^\sigma(e)=\frac34+\sigma\frac{\sqrt{33}}4e+O(e^2),
\qquad \sigma\in\{+1,-1\}.
\]
The simple crossings above imply continuity of the unique branches $\rho^\pm(e)$. Since they never meet for $0<e<1$, their order near $e=0$ persists throughout this interval. Hence
\[
\beta_s(e)=\rho^-(e),\qquad
\beta_m(e)=\rho^+(e),\qquad 0<e<1.
\]
This completes the proof of Theorem~\ref{thm:main}.

\section{Concluding remarks}
Both strict separation results originate in the reflection dependence of the lowest half-integer Fourier block. The sign of the resulting perturbation determines the mechanism: in the Euler family, a definite connecting crossing form permits a spectral-flow comparison; in the Lagrangian family, the indefinite first-block difference is handled by scalar control of the common tail and an effective-block estimate. The kernel transformations preserve degeneracy in both arguments. The Euler proof uses the global periodic separators and sector counts, while the Lagrangian proof bounds the tail correction to the first-block diagonal difference.

The conclusions hold for every $0<e<1$. Estimate~\eqref{eq:effectivegap} controls the diagonal difference of the effective Lagrangian block; quantitative bounds for the parameter gap $\beta_m(e)-\beta_s(e)$ would involve further parameter estimates. The singular limit $e\to1$ is a natural direction for further study.

\FloatBarrier
\section*{Acknowledgements}
This work was supported by the Taishan Scholars Climbing Program of Shandong (TSPD20240802).

Generative AI tools were used to assist with coding, computational analysis, and manuscript preparation. The authors directed and verified the work and take full responsibility for the content of the manuscript.

% Bibliography embedded for direct TeXShop compilation.

\end{document}

%% file: figures/euler-proof-panels.tex
% Requires tikz and \usetikzlibrary{arrows.meta}.
% Proof schematic only: the local zero curve is not a computed branch.
\begingroup
\definecolor{eulerProofFlow}{RGB}{35,65,91}
\definecolor{eulerProofLocal}{RGB}{22,122,129}
\begin{tikzpicture}[
  x=1cm,y=1cm,font=\footnotesize,
  >={Latex[length=1.7mm,width=1.2mm]},
  flow/.style={draw=eulerProofFlow,line width=.65pt,->},
  local/.style={draw=eulerProofLocal,line width=.8pt},
  guide/.style={draw=gray!65,line width=.35pt,densely dashed},
  axis/.style={draw=gray!70,line width=.4pt,->}
]
\path[use as bounding box] (0,-.45) rectangle (15.7,5.1);

% Panel (a): the arrows specify traversal of the rectangle boundary.
\node[font=\small] at (3.75,4.65) {(a) Spectral-flow rectangle};
\draw[flow] (1.8,1.25)--(5.45,1.25)
  node[midway,below=5pt] {$+N_-(\beta)$};
\draw[flow] (5.45,3.65)--(1.8,3.65)
  node[midway,above=5pt] {$-N_+(\beta)$};
\draw[flow,draw=eulerProofLocal] (5.45,1.25)--(5.45,3.65)
  node[midway,right=6pt] {$+\Phi(\beta)\ge0$};
\draw[flow,draw=eulerProofLocal] (1.8,3.65)--(1.8,1.25)
  node[midway,left=6pt,align=center] {$-\Phi(a)$\\$=0$};
\node[left] at (1.65,1.25) {$t=-1$};
\node[left] at (1.65,3.65) {$t=1$};
\draw[axis] (1.45,.55)--(6.0,.55);
\node[below=3pt] at (1.8,.55) {$a$};
\node[below=3pt] at (5.45,.55) {$\beta$};
\node at (3.75,-.25) {$\Phi(\beta)=N_+(\beta)-N_-(\beta)$};

% Panel (b): only the neighbourhood of the lower-sector endpoint is shown.
\node[font=\small] at (11.8,4.65) {(b) Local zero curve};
\draw[axis] (8.9,.55)--(15.25,.55);
\draw[axis] (9.35,.55)--(9.35,4.05) node[above] {$t$};
\node[right=2pt] at (15.25,.55) {$\beta$};
\draw[guide] (9.35,1.15)--(15.05,1.15);
\node[left=3pt] at (9.35,1.15) {$-1$};
\draw[guide] (12.3,.55)--(12.3,2.65);
\draw[guide] (13.8,.55)--(13.8,1.15);
\node[below=3pt] at (12.3,.55) {$\beta_*$};
\node[below=3pt] at (13.8,.55) {$\rho^-$};
\draw[local] plot[domain=11.7:13.8,samples=30]
  (\x,{1.15+.85*(13.8-\x)+.1*(13.8-\x)^2});
\draw[local,->] plot[domain=13.3:12.7,samples=12]
  (\x,{1.15+.85*(13.8-\x)+.1*(13.8-\x)^2});
\fill[eulerProofLocal] (12.3,2.65) circle (1.5pt);
\fill[eulerProofFlow] (13.8,1.15) circle (1.5pt);
\node[above right=3pt] at (12.3,2.65) {$(\beta_*,t_*)$};
\node[above right=3pt] at (13.8,1.15) {$(\rho^-,-1)$};
\node[anchor=east] at (11.7,2.65) {$\lambda(\beta,t)=0$};
\node[anchor=west] at (11.75,3.65) {$d'(-1)<0$};
\node at (11.8,-.25) {$t_*>-1,\qquad \beta_*<\rho^-$};
\end{tikzpicture}
\endgroup

%% file: main.bbl
\begin{thebibliography}{99}
\small
\raggedright
\setlength{\itemsep}{4pt plus 1pt minus 0.5pt}

\bibitem{AtiyahPatodiSinger1976}
M.~F. Atiyah, V.~K. Patodi and I.~M. Singer,
\href{https://doi.org/10.1017/S0305004100052105}{Spectral asymmetry and {Riemannian} geometry. {III}},
\emph{Math. Proc. Cambridge Philos. Soc.} \textbf{79} (1976), 71--99.

\bibitem{FitzpatrickPejsachowiczRecht1999}
P.~M. Fitzpatrick, J. Pejsachowicz and L. Recht,
\href{https://doi.org/10.1006/jfan.1998.3366}{Spectral flow and bifurcation of critical points of strongly-indefinite functionals. Part {I}. General theory},
\emph{J. Funct. Anal.} \textbf{162} (1999), 52--95.

\bibitem{HuLiuOuQiaoSalomao2026}
X. Hu, L. Liu, Y. Ou, Z. Qiao and P.~A.~S. Salom{\~a}o,
\href{https://doi.org/10.2140/gt.2026.30.3513}{{ECH} constraints and twist dynamics in the spatial isosceles three-body problem},
\emph{Geom. Topol.} \textbf{30} (2026), no.~9, 3513--3574.

\bibitem{HuLiuOuSalomaoYu2025}
X. Hu, L. Liu, Y. Ou, P.~A.~S. Salom{\~a}o and G. Yu,
\href{https://doi.org/10.4171/JEMS/1678}{A symplectic dynamics approach to the spatial isosceles three-body problem},
\emph{J. Eur. Math. Soc.} (2025),
\href{https://doi.org/10.4171/JEMS/1678}{doi:10.4171/JEMS/1678}.

\bibitem{HuLongSun2014}
X. Hu, Y. Long and S. Sun,
\href{https://doi.org/10.1007/s00205-014-0749-6}{Linear stability of elliptic {Lagrangian} solutions of the planar three-body problem via index theory},
\emph{Arch. Ration. Mech. Anal.} \textbf{213} (2014), 993--1045.

\bibitem{HuOu2016}
X. Hu and Y. Ou,
Collision index and stability of elliptic relative equilibria in planar {$n$}-body problem,
\emph{Comm. Math. Phys.} \textbf{348} (2016), 803--845.

\bibitem{HuSun2009Sym}
X. Hu and S. Sun,
Index and stability of symmetric periodic orbits in {Hamiltonian} systems with application to figure-eight orbit,
\emph{Comm. Math. Phys.} \textbf{290} (2009), 737--777.

\bibitem{HuSun2009RE}
X. Hu and S. Sun,
Stability of relative equilibria and {Morse} index of central configurations,
\emph{C. R. Math. Acad. Sci. Paris} \textbf{347} (2009), 1309--1312.

\bibitem{HuSun2010}
X. Hu and S. Sun,
{Morse} index and stability of elliptic {Lagrangian} solutions in the planar three-body problem,
\emph{Adv. Math.} \textbf{223} (2010), 98--119.

\bibitem{Long2002}
Y. Long,
\emph{Index Theory for Symplectic Paths with Applications},
Progr. Math., vol.~207, Birkh{\"a}user, Basel, 2002.

\bibitem{MSS2006b}
R. Mart{\'i}nez, A. Sam{\`a} and C. Sim{\'o},
Analysis of the stability of a family of singular-limit linear periodic systems in {$\mathbb{R}^4$}. Applications,
\emph{J. Differential Equations} \textbf{226} (2006), 652--686.

\bibitem{MSS2006a}
R. Mart{\'i}nez, A. Sam{\`a} and C. Sim{\'o},
Stability diagram for {4D} linear periodic systems with applications to homographic solutions,
\emph{J. Differential Equations} \textbf{226} (2006), 619--651.

\bibitem{MeyerSchmidt2005}
K.~R. Meyer and D.~S. Schmidt,
Elliptic relative equilibria in the {$N$}-body problem,
\emph{J. Differential Equations} \textbf{214} (2005), 256--298.

\bibitem{Phillips1996}
J. Phillips,
\href{https://doi.org/10.4153/CMB-1996-054-4}{Self-adjoint {Fredholm} operators and spectral flow},
\emph{Canad. Math. Bull.} \textbf{39} (1996), 460--467.

\bibitem{RobbinSalamon1993}
J. Robbin and D. Salamon,
\href{https://doi.org/10.1016/0040-9383(93)90052-W}{The {Maslov} index for paths},
\emph{Topology} \textbf{32} (1993), 827--844.

\bibitem{RobbinSalamon1995}
J. Robbin and D. Salamon,
\href{https://doi.org/10.1112/blms/27.1.1}{The spectral flow and the {Maslov} index},
\emph{Bull. Lond. Math. Soc.} \textbf{27} (1995), 1--33.

\bibitem{Roberts2002}
G.~E. Roberts,
Linear stability of the elliptic {Lagrangian} triangle solutions in the three-body problem,
\emph{J. Differential Equations} \textbf{182} (2002), 191--218.

\bibitem{ZhouLong2017}
Q. Zhou and Y. Long,
\href{https://doi.org/10.1007/s00205-017-1154-8}{Maslov-type indices and linear stability of elliptic {Euler} solutions of the three-body problem},
\emph{Arch. Ration. Mech. Anal.} \textbf{226} (2017), 1249--1301.

\end{thebibliography}
